\documentclass[12pt]{article}
\usepackage{a4wide}
\usepackage{amsmath}
\usepackage{amssymb}
\usepackage{amsthm}
\usepackage{mathtools}
\usepackage{mdwlist}
\usepackage{color}
\usepackage{graphicx}

\newtheorem{lemma}{Lemma}

\newtheorem{remark}{Remark}
\newtheorem{definition}{Definition}
\newtheorem{proposition}{Proposition}

\newcommand{\eps}{\varepsilon}

\renewcommand{\rho}{\varrho}
\renewcommand{\phi}{\varphi}
\renewcommand{\theta}{\vartheta}

\newcommand{\abs}[1]{\left|#1\right|}
\newcommand{\norm}[1]{\left\|#1\right\|}

\newcommand{\ms}{\hspace*{-0.1em}}
\newcommand{\en}[1]{\left|\ms\left|\ms\left|#1\right|\ms\right|\ms\right|}

\begin{document}

\title{An $hp$ Weak Galerkin FEM for singularly perturbed problems}
\author{Torsten Lin{\ss }\thanks{%
LG Numerische Mathematik, Fakult\"at f\"ur Mathematik und Informatik,
FernUniversit\"at in Hagen, D-58084 Hagen, Germany (\texttt{%
torsten.linss@fernuni-hagen.de}).} \hspace{0.025cm} and Christos Xenophontos%
\thanks{%
Department of Mathematics \& Statistics, University of Cyprus, PO BOX 20537,
Nicosia 1678, Cyprus (\texttt{xenophontos@ucy.ac.cy}).} }
\maketitle

\begin{abstract}
We present the analysis for an $hp$ weak Galerkin-FEM for singularly
perturbed reaction-convection-diffusion problems in one-dimension. 
Under the analyticity of the data
assumption, we establish robust exponential convergence, when the error is
measured in the energy norm, as the degree $p$ of the approximating
polynomials is increased. The {\emph{Spectral Boundary Layer}} mesh is used,
which is the minimal (layer adapted) mesh for such problems. Numerical
examples illustrating the theory are also presented.
\end{abstract}

\section{Introduction}

In recent years there has been a growing interest in {\emph{hybrid}}
methods, which involve modified variational formulations and/or
discretizations (see \cite{CEP} and the references therein). One such method
is the so-called {\emph{Weak Galerkin}} (WG) method \cite{WY}, which is closely related 
to the Discontinuous Galerkin (DG) method \cite{DiPE}. In \cite{ZX}, 
singularly perturbed convection-diffusion problems were considered and a WG $h$ 
version discretization was presented and analyzed (see also \cite{XZZ}).
In the present article we extend the method and analysis in \cite{ZX} as follows: we consider
a more general reaction-convection-diffusion problem (with two small parameters) and we use
the $p$ version in conjunction with the {\emph{Spectral Boundary Layer}} mesh \cite{MXO}. 
We show that the method converges exponentially and independently of the singular perturbation parameters, 
as the degree $p$ of the approximating polynomials is increased and the error
is measured in the energy norm.

The rest of the article is organized as follows: in Section \ref{sec:mesh} we describe the 
discrete setting with the mesh and spaces utilized, in Section \ref{sec:model} we present
the model problem and its regularity, and in Section \ref{sec:WG} we give the 
{\emph{Weak Galerkin}} formulation and related results. The approximation for
the three different cases (reaction-diffusion, convection-diffusion and reaction-convection-diffusion)
is presented in Section \ref{sec:approx}, and in Section \ref{NR} we illustrate our theoretical
findings through numerical computations.

Throughout the article, we will utilize standard Sobolev spaces:
With $I \subset \mathbb{R},\ $a domain with boundary $\partial I $ and measure 
$\abs{I}$, we will denote by $C^{k}(I )$ the space of continuous functions on $I $ with
continuous derivatives up to order $k$. We will use the Sobolev spaces $W^{k,m}(I )$ of 
functions on $I$ with $0,1,2,...,k$ generalized derivatives in $L^{m}\left(I \right) $, equipped 
with the norm and seminorm $\left\Vert \cdot \right\Vert _{k,m,I }$ and
 $\abs{\cdot} _{k,m,I}\,$, respectively. When $m=2$, we will write 
$H^{k}\left( I \right) $ instead of $W^{k,2}\left(I \right) $, and
for the norm and seminorm, we will write $\left\Vert \cdot \right\Vert
_{k,I }$ and $\abs{\cdot}_{k,I }\,$,
respectively. The usual $L^{2}(I )$ inner product will be denoted by 
$\left ( \cdot ,\cdot \right ) _{I }$, with the subscript
omitted when there is no confusion. We will also use the space 
\[
H_{0}^{1}\left(I \right) =\left\{ u\in H^{1}\left( \Omega \right)
:\left. u\right\vert _{\partial \Omega }=0\right\} .
\]
The norm of the space $L^{\infty }(I)$ of essentially bounded
functions is denoted by $\Vert \cdot \Vert _{\infty ,I }$. Finally, the
notation \textquotedblleft $a\lesssim b$\textquotedblright\ means
\textquotedblleft $a\leq Cb$\textquotedblright\ with $C$ being a generic
positive constant, independent of any discretization or singular
perturbation parameters.


\section{Mesh and discrete spaces}
\label{sec:mesh}
Let $\omega$ be an arbitrary mesh on $I=(0,1)$, i.e. 
\[
\omega \colon 0 = x_0 < x_1 < \cdots < x_N = 1, \ \ I_j \coloneqq \bigl(%
x_{j-1}, x_j\bigr) , \ \ h_j \coloneqq x_{j} - x_{j-1}.
\]
We say that $v:=\{v_0,v_b\}$ is a {\emph{weak function}} on $I_j$, if $v_0
\in L^2(I_j)$ and $v_b \in L^{\infty}(\partial I_j)$, where $\partial I_j =
\{x_{j-1},x_j\}$ denotes the endpoints (or the boundary) of $I_j$. In
particular, 
\begin{align*}
v = \left\{ 
\begin{array}{ccc}
v_0 \;,\; & {\text{in}} \; & I_J \\ 
v_b \; , \; & {\text{on}} \; & \partial I_j%
\end{array}
\right.
\end{align*}
with $v_0$ giving the values in $(x_{j-1},x_j)$ and $v_b$ giving the values
at the endpoints. The space of weak functions on $I_j$ is denoted by 
\begin{equation*}
\mathcal{M}(I_j) = \{ v=\{v_0,v_b\}: v_0 \in L^2(I_j), v_b \in
L^{\infty}(\partial I_j) \}.
\end{equation*}
We also define the following function spaces which will be utilized in the
sequel: 
\begin{align*}
\mathcal{M}^\omega & \coloneqq \left\{\bigl(v_0, v_b\bigr) \in L_2(I) \times 
\mathbb{R}^{N+1} \colon v_0 \vert_{I_{j}} \in H^1(I_{j}) , \ j=1,\ldots,N; \
v_b = (v_{b,0}, \ldots , v_{b,N}) \right\}, \\
\mathcal{M}^\omega_{0} & \coloneqq \left\{(v_0, v_b) \in \mathcal{M}%
^{\omega} \colon v_{b,0} = v_{b,N} = 0 \right\}.
\end{align*}
Note that for functions $v\in H^1(I)\cap \mathcal{M}^\omega$ we have $%
v_{b,j} = v_0 (x_{j})$. Moreover, for functions $v=(v_0,v_b)\in\mathcal{M}%
^\omega$ we set $v_{0,j} \coloneqq v_0\vert_{I_{j}}$. With $\mathbb{P}_p$
the space of polynomials of degree $\leq p$, we define the FE-spaces: 
\begin{align*}
V^\omega_p & \coloneqq \left\{(v_0,v_b) \in \mathcal{M}^{\omega} \colon v_0
\mid_{I_{j}} \in \mathbb{P}_{p}, \ j = 1,\ldots ,N \right\} \subset \mathcal{%
M}^{\omega}, \\
V^\omega_{p,0} & \coloneqq \left\{(v_0,v_b) \in V^\omega_p \colon v_{b,0} =
v_{b,N} = 0 \right\} = V^\omega_p \cap \mathcal{M}^\omega_0, \\
W^\omega_p & \coloneqq \left\{ v_{0} \in L_{2} (0,1) \colon v_0 \mid_{I_{j}}
\in \mathbb{P}_p, \ j = 1,\ldots ,N \right\}.
\end{align*}

We next present definitions appearing in \cite{ZX}, that describe the
notions of {\emph{weak derivative}} and {\emph{weak convection derivative}},
respectively.

\begin{definition}
Let $p\in\mathbb{N}_0$. The weak derivative of a function $v \in \mathcal{M}%
^{\omega}$, is the unique polynomial $D_p v \in W^\omega_p$, satisfying %
%
%
\[
\int_{I_{j}} D_p v q = - \int_{I_{j}} v_0 q' + \langle v_b , q \cdot n
\rangle_{\partial I_{j}} \,, \ \quad \forall \ q \in \mathbb{P}_{p} \,, \ j
= 1,\ldots ,N \,,
\]
where $\langle v_b,q \cdot n \rangle_{\partial I_{j}} \coloneqq v_{b,j} q
(x_{j}) - v_{b,j-1} q(x_{j-1})$.
\end{definition}

\begin{definition}
Let $p\in\mathbb{N}_0$. The weak convection derivative of a function $v \in 
\mathcal{M}^{\omega}$, is the unique polynomial $D^c_p v \in W^\omega_p$,
satisfying %
%
%
\[
\int_{I_{j}} D^c_p v q = - \int_{I_{j}} v_0 (bq)^{\prime }+ \langle v_b ,
(bq) \cdot n \rangle_{\partial I_{j}} \,, \ \quad \forall \ q \in \mathbb{P}%
_{p} \,, \ j = 1,\ldots ,N \,,
\]
where $\langle v_b,(bq) \cdot n \rangle_{\partial I_{j}} \coloneqq v_{b,j}
(bq) (x_{j}) - v_{b,j-1} (bq)(x_{j-1})$, and $b$ is any sufficiently smooth function.
\end{definition}

The above will appear in the {\emph{weak Galerkin discretization}}, which
amounts to using weak derivatives in the variational formulation (see Section \ref{sec:WG} ahead).


\section{The model problem}
\label{sec:model}
We consider the following problem: find $u$ such that 
\begin{eqnarray}
- \varepsilon_1 u^{\prime \prime }(x) + \varepsilon_2 b(x) u^{\prime }(x) +
r(x)u(x) &= f(x) \;,\; x\in I=(0,1)  \label{de} \\
u(0)=u(1)=0  \label{bc}
\end{eqnarray}
where $0<\varepsilon_1, \varepsilon_2 \leq 1$ are given parameters that can
approach zero and the functions $b, r,f$ are given and sufficiently smooth;
in particular, we assume that they are {\emph{analytic}}
functions satisfying, for some positive constants $\gamma_b, \gamma
_{r},\gamma _{f}$, independent of $\varepsilon$, 
\begin{equation}
\left\Vert b^{(n)}\right\Vert _{L^{\infty}(I)}\lesssim
n!\gamma_{b}^{n}, \left\Vert r^{(n)}\right\Vert
_{L^{\infty}(I)}\lesssim n!\gamma_{r}^{n}, \left\Vert
f^{(n)}\right\Vert _{L^{\infty}(I)}\lesssim n!\gamma_{f}^{n}, \;\forall
\;n\in \mathbb{N}_{0}.  \label{analytic}
\end{equation}
Moreover, we assume that $b(x)>0, r(x) \ge 0 \; \forall \; x \in [0,1]$,
and that there exists a positive constant $\gamma$ such that 
\begin{gather}
r - \frac{\eps_2 b^{\prime }}{2} \ge \gamma > 0 \quad \text{on} \ [0,1].
\end{gather}
(These are the usual assumptions made in order for the variational problem to have 
a unique solution, see e.g. \cite{LR})

More details arise if one studies the structure of the solution to (\ref{de}%
), which depends on the roots of the characteristic equation associated with
the differential operator. For this reason, we let $\lambda_{0}(x),\lambda
_{1}(x)$ be the solutions of the characteristic equation and set 
\begin{equation}
\mu _{0}=-\underset{x\in \lbrack 0,1]}{\max }\lambda _{0}(x)\;,\;\mu _{1}=%
\underset{x\in \lbrack 0,1]}{\min }\lambda _{1}(x),  \label{mu}
\end{equation}
or equivalently, 
\begin{equation*}
\mu _{0,1}=\underset{x\in \lbrack 0,1]}{\min }\frac{\mp \varepsilon _{2}b(x)+%
\sqrt{\varepsilon _{2}^{2}b^{2}(x)+4\varepsilon _{1}r(x)}}{2\varepsilon _{1}}%
.
\end{equation*}
The following hold true \cite{Linss,RU}:%
\begin{equation}
\left. 
\begin{array}{c}
1\ll\mu _{0}\leq \mu _{1\;\;},\;\frac{\varepsilon _{2}}{\varepsilon
_{2}+\varepsilon _{1}^{1/2}}\lesssim \varepsilon _{2}\mu _{0}\lesssim
1\;,\;\varepsilon _{1}^{1/2}\mu _{0}\lesssim 1 \\ 
\max \{\mu _{0}^{-1},\varepsilon _{1}\mu _{1}\}\lesssim \varepsilon
_{1}+\varepsilon _{2}^{1/2}\;,\;\varepsilon _{2}\lesssim \varepsilon _{1}\mu
_{1} \\ 
\text{for }\varepsilon _{2}^{2}\geq \varepsilon _{1}:\;\varepsilon
_{1}^{-1/2}\lesssim \mu _{1}\lesssim \varepsilon _{1}^{-1} \\ 
\text{for }\varepsilon _{2}^{2}\leq \varepsilon _{1}:\;\varepsilon
_{1}^{-1/2}\lesssim \mu _{1}\lesssim \varepsilon _{1}^{-1/2}%
\end{array}
\right\} .  \label{mu_a}
\end{equation}
The values of $\mu _{0},\mu _{1}$ determine the strength of the boundary
layers and since $\left\vert \lambda _{0}(x)\right\vert <\left\vert \lambda
_{1}(x)\right\vert $ the layer at $x=1$ is stronger than the layer at $x=0$.
Essentially, there are three regimes \cite{Linss}, \cite{LR}, as summarized in Table 1.

\begin{table}[h]
\begin{center}
\begin{tabular}{||cccc||}
\hline
&  & $\mu _{0}$ & $\mu _{1}$ \\[0.5ex] \hline\hline
convection-diffusion & $\varepsilon _{1}\ll\varepsilon _{2}=1$ & $1$ & $%
1/\varepsilon _{1}$ \\ \hline
convection-reaction-diffusion & $\varepsilon _{1}\ll\varepsilon _{2}^{2}\ll1$
& $1/\varepsilon _{2}$ & $\varepsilon _{2}/\varepsilon _{1}$ \\ \hline
reaction-diffusion & $1\gg\varepsilon _{1}\gg\varepsilon _{2}^{2}$ & $1/\sqrt{%
\varepsilon _{1}}$ & $1/\sqrt{\varepsilon _{1}}$ \\[1ex] \hline
\end{tabular}%
\end{center}
\caption{\label{table1}
Different regimes based on the relationship between $\protect%
\varepsilon _{1}$ and $\protect\varepsilon _{2}$.}
\end{table}

In \cite{SXArXiv} it was shown that
\begin{equation}
u=u_{S}+u^{\pm}_{BL}+u_{R},  \label{decomposition}
\end{equation}%
where $u_{S}$ denotes the smooth part, $u^{\pm}_{BL}$ denote the boundary layers
at the two endpoints, and $u_{R}$ denotes the remainder. Estimates on each component may 
be established, which are explicit in the singular perturbation parameters and the order of 
differentiation. More details will be given in Section \ref{sec:reg} below.

The (standard Galerkin) variational formulation of (\ref{de})--(\ref{bc}) reads: Find $u \in
H^1_0(I)$ such that 
\begin{gather*}
{\mathcal{A}}(u,v) \coloneqq \varepsilon_1(u^{\prime },v^{\prime })_I + \varepsilon_2(b u^{\prime },v)_I
+(ru,v)_I = (f,v)_I \quad \forall \ v \in H^1_0(I) .
\end{gather*}

\subsection{Regularity}
\label{sec:reg}
As seen in Table 1, we have three different types of problems, based on the relationship between
$\eps_1, \eps_2$. In this section, we present analytic regularity results for each case, as derived in 
\cite{melenk97}, \cite{MS}, \cite{SXArXiv}, under the assumption (\ref{analytic}) of analytic input data.

\subsubsection{Reaction-convection-diffusion}
In this case $\eps_1 \ll \eps_2^2$, and we have layers at each endpoint of different widths: $O(\eps_2)$ at the left and $O(\eps_2 / \eps_1)$ 
at the right. In \cite{SXArXiv} it was shown that $u$ saitisfies (\ref{decomposition}) and there exist positive constants  $K_{1}$, $\tilde{K}_S$, $\tilde{K}$, $\bar{K}$, $\tilde{\delta}$, independent of $\varepsilon _{1},\varepsilon _{2},$ such that $\forall \; n \in \mathbb{N}$ and $x \in [0,1]$, there holds
\begin{eqnarray}
\Vert u^{(n)} \Vert_{L^{\infty}(I)} &\lesssim& {K}_1^n \max\{ n, \eps_2/ \eps_1, \eps_2^{-1} \}^{n}  \; ,  \label{Classical3} \\
\left\Vert u_{S}^{(n)}\right\Vert _{\infty ,I} &\lesssim&  n! \tilde{K}_{S}^{n}\; ,   \label{uM1bound} \\
\left\vert \left(u_{BL}^{-}\right) ^{(n)}(x)\right\vert & \lesssim & 
\tilde{K}^{n}\varepsilon _{2}^{-n}e^{-dist(x,\partial I)/\varepsilon_{2}}\; ,  \label{uBL1abound} \\
\left\vert \left( u_{BL}^{+}\right) ^{(n)}(x)\right\vert &\lesssim& \bar{K}^{n}\left( \frac{\varepsilon _{1}}{\varepsilon _{2}}\right)
^{-n}e^{-dist(x,\partial I)\varepsilon _{2}/\varepsilon _{1}}\; ,   \label{uBL1bbound} \\
\left\Vert u_{R} \right\Vert _{\infty ,\partial I}+\left\Vert  u_{R} \right\Vert _{0,I}+\varepsilon _{1}^{1/2}\left\Vert u_{R}^{\prime
}\right\Vert _{0,I} &\lesssim& \max \{e^{-\tilde{\delta} \varepsilon _{2}/\varepsilon
_{1}},e^{-\tilde{\delta} /\varepsilon _{2}}\}.  \label{rM1bound}
\end{eqnarray}
%

\subsubsection{Reaction-diffusion}

In this case $\eps_1 \gg \eps_2^2$, and we have two layers of width $O(\eps_1^{1/2})$ at the endpoints.
It was shown in \cite{melenk97} that the solution to pproblem (\ref{de})--(\ref{bc}) admits the decomposition
(\ref{decomposition}), and there exist positive constants $K, K_S, K_{\pm}, \delta$ such that
 $\forall \; n \in \mathbb{N}$ and $x \in [0,1]$, there holds
\begin{eqnarray}
\Vert u^{(n)} \Vert_{L^{\infty}(I)} &\lesssim& K^n \max\{ n, \eps_1^{-1/2}\}^n  \; , \label{Classical} \\
\Vert u_S^{(n)} \Vert_{L^{\infty}(I)} &\lesssim& n! K_S^n \label{S_bound} \; ,  \\
\vert (u^{\pm}_{BL})^{(n)} (x)\vert &\lesssim& K_{\pm}^n \eps_1^{-n/2} e^{-\text{dist}(x,\partial I)/\eps_1^{1/2}}  \; ,  \label{BL_bound} \\
\vert u_{R}(\pm 1) \vert + \eps_1 \vert u_{R} \vert_{H^1(I)} + \Vert u_{R} \Vert_{L^2(I)} &\lesssim& e^{-\delta / \eps_1^{1/2}}. \label{R_bound1}
\end{eqnarray}

\subsubsection{Convection-diffusion}
In this case $\eps_2 = 1$, and we have a boundary layer only at the outflow, i.e. $x=1$, of width $O(\eps_1)$. In \cite{MS} it was shown that 
 there exist positive constants $\hat{K}, \hat{K}_S, \hat{K}_{BL}, \hat{\delta}$ such that $u$ saitisfies (\ref{decomposition}) and
 $\forall \; n \in \mathbb{N}$ and $x \in [0,1]$, there holds
\begin{eqnarray}
\Vert u^{(n)} \Vert_{L^{\infty}(I)} &\lesssim& \hat{K}^n \max\{ n, \eps_1^{-1}\}^n  \; , \label{Classical2} \\
\Vert u_S^{(n)} \Vert_{L^{\infty}(I)} &\lesssim&  n! \hat{K}_S^n \; , \label{S_bound2} \\
\vert (u_{BL})^{(n)} (x)\vert &\lesssim& \hat{K}_{BL}^n \eps_1^{-n} e^{-\text{dist}(x,\partial I)/\eps_1}  \; ,\label{BL_bound2} \\
\vert u_{R}(\pm 1) \vert + \eps_1 \vert u_{R} \vert_{H^1(I)} + \Vert u_{R} \Vert_{L^2(I)} &\lesssim& e^{-\hat{\delta} / \eps_1}.
\end{eqnarray}

The above regularity results will be used in the design of the mesh as well as in the numerical analysis of the method.

\subsection{Weak Galerkin formulation}
\label{sec:WG}
The \emph{Weak Galerkin} (WG) discretization reads as follows: Given $p\in \mathbb{N}$, 
find $u_{p}=(u_{0},u_{b})\in V_{p,0}^{\omega }$ such that, 
\begin{equation}
\label{Auv}
{\mathcal{A}}_{p}(u_{p},v) = (f,v_{0})_{I} \;\; \forall \ v=(v_{0},v_{b})\in V_{p,0}^{\omega },
\end{equation}
where 
\begin{equation}
\label{Ap}
{\mathcal{A}}_{p}(u_{p},v)\coloneqq\varepsilon
_{1}(D_{p-1}u_{p},D_{p-1}v)_{I}+\varepsilon
_{2}(D_{p}^{c}u_{p},v_{0})_I +(ru_{0},v_{0})_{I}+S(u,v)+S_{c}(u,v)
\end{equation}
with
\begin{equation} \label{Suv}
S(u,v)\coloneqq\sum_{j=1}^{N}\sigma _{j}\langle
u_{0}-u_{b},v_{0}-v_{b}\rangle _{\partial I_{j}},\quad \langle \phi ,\psi
\rangle _{\partial I_{j}}\coloneqq(\phi \psi )(x_{j}^{-})+(\phi \psi
)(x_{j-1}^{+}),
\end{equation}
and
\begin{equation} \label{Scuv}
S_{c}(u,v)\coloneqq\sum_{j=1}^{N}\langle
u_{0}-u_{b},(v_{0}-v_{b}) \eps_2 b\cdot n_{j}\rangle _{\partial _{+}I_{j}},\quad
\langle \phi ,\psi \cdot n_{j}\rangle _{\partial _{+}I_{j}}\coloneqq(\phi
\psi )(x_{j}^{-})
\end{equation}
where the constants $\sigma _{j}$, called {\emph{penalty parameters}}, are chosen 
as $\eps_1 p^2 / h_j$ (as was done in \cite{HSS} for $hp$ DG-FEM).

We also introduce the norms $\forall \; v=(v_{0},v_{b})\in \mathcal{M}^{\omega }$:
\begin{eqnarray}\en{v}_{p, I}^{2 } &\coloneqq& \varepsilon _{1}\left\Vert
D_{p-1}v\right\Vert _{0,I}^{2}+\left\Vert v_{0}\right\Vert
_{0,I}^{2}+S(v,v)+S_{c}(v,v)+\left\vert v\right\vert _{J}^{2}, \label{p_norm} \\
\en{v} ^{2} &\coloneqq& \varepsilon _{1}\left\Vert v_{0}^{\prime }\right\Vert
_{0,\omega }^{2}+\left\Vert v_{0}\right\Vert
_{0,I}^{2}+S(v,v)+S_{c}(v,v)+\left\vert v\right\vert _{J}^{2}, \label{tr_norm}
\end{eqnarray}
where
\[
\left\Vert v_{0}^{\prime }\right\Vert _{0,\omega }^{2} \coloneqq
\sum_{j=1}^{N}\left\Vert v_{0,j}^{\prime }\right\Vert _{0,I_{j}}^{2}
\] 
and
\begin{equation*}
\left\vert v\right\vert _{J}^{2}\coloneqq%
\sum_{j=1}^{N}w_{j} \eps_2 b(x^-_{j})(v_{0}-v_{b})(x_{j}^{-})^{2}\,,\quad w_{j}%
\coloneqq%
\begin{cases}
1 & j=1,\dots ,N-1, \\ 
1/2 & j=N%
\end{cases}%
\end{equation*}
The subscript $I$ will be omitted when there is no confusion.

The norms  (\ref{p_norm}), (\ref{tr_norm}) are equivalent, under a certain condition as
stated in the following.

\begin{lemma}
Under the assumption that there exists a constant $C_{\sigma }>0$ such that 
\begin{equation} \label{sigma_j}
\varepsilon _{1} p^2 h_{j}^{-1}/\sigma _{j}\leq C_{\sigma }\;,\;j=1,\ldots ,N,
\end{equation}
the norms $\left\vert \hspace*{-0.1em}\left\vert \hspace*{-0.1em}\left\vert
\cdot \right\vert \hspace*{-0.1em}\right\vert \hspace*{-0.1em}\right\vert
_{p}$ and $\left\vert \hspace*{-0.1em}\left\vert \hspace*{-0.1em}\left\vert
\cdot \right\vert \hspace*{-0.1em}\right\vert \hspace*{-0.1em}\right\vert $
are equivalent on $V_{p}^{\omega }$, i.e. there exist constants $C_{lb}$ and 
$C_{ub}$ such that 
\begin{equation*}
C_{lb}\left\vert \hspace*{-0.1em}\left\vert \hspace*{-0.1em}\left\vert
v\right\vert \hspace*{-0.1em}\right\vert \hspace*{-0.1em}\right\vert \leq
\left\vert \hspace*{-0.1em}\left\vert \hspace*{-0.1em}\left\vert
v\right\vert \hspace*{-0.1em}\right\vert \hspace*{-0.1em}\right\vert
_{p}\leq C_{ub}\left\vert \hspace*{-0.1em}\left\vert \hspace*{-0.1em}%
\left\vert v\right\vert \hspace*{-0.1em}\right\vert \hspace*{-0.1em}%
\right\vert \quad \forall v=(v_{0},v_{b})\in V_{p}^{\omega }\ .
\end{equation*}
The constants are independent of the mesh $\omega $ and the polynomial
degree $p$.
\end{lemma}
\begin{proof}
This was shown in \cite{ZX} for the convection-diffusion case. Inspecting the proof, we
see that the result holds in the more general case as well, with only minor modifications.
\end{proof}

The next lemma gives the coercivity of the bilinear form (\ref{Auv}).

\begin{lemma}
  The bilinear form ${\cal{A}}_p$ defined by (\ref{Auv}) is coercive with respect to
  both norms $\en{ \cdot }$ and $\en{ \cdot}_p$:
  \begin{alignat*}{2}
    \mathcal{A}_p (v,v) & \geq \en{v}^2  & \forall & \; v \in {\cal{M}}_0^{\omega}, \\
    \mathcal{A}_p (v,v) & \geq C_{lb} \en{v}^2_p \ \ & \forall & \; v \in {\cal{M}}_0^{\omega}.
\end{alignat*}
\end{lemma}
\begin{proof}
As with the previous lemma, the result was established in \cite{ZX} for the convection-diffusion case. The extension to this 
case is straight forward.
\end{proof}

\subsection{Interpolation results}

We next describe the interpolant we will use:

\begin{definition}
\label{interpolant} We define the linear operator ${\mathcal{I}}\colon
H^{p+1}(0,1) \to V^\omega_p \colon v \mapsto {\mathcal{I}} v$, by the
conditions 
\begin{gather*}
\bigl({\mathcal{I}} v -v\bigr)(x_j) = 0 \,, \ \int_{I_j} ({\mathcal{I}}
v-v)^{\prime }q^{\prime }= 0 \quad \forall \ q \in \mathbb{P}_p, \ j =
1,\ldots,N\,.
\end{gather*}
\end{definition}

\begin{remark}
The above interpolant is identical to the one obtained via \textsc{Legendre}%
-expansions, see \cite{Schwab98}. Consequently, the following lemma is shown
there (see \cite[Theorem 3.15]{Schwab98}), and will be the main tool in our analysis.
\end{remark}

\begin{lemma}\label{prop:schwab}
Let $\hat{I}=(a,b)$ and assume  $y \in H^k(\hat{I})$ for some $k \ge 1$,
with ${\mathcal{I}} y$ its interpolant given by Definition \ref{interpolant}. Then, 
\begin{gather*}
\left|y-{\mathcal{I}} y\right|_{1,\hat{I}} + k \left\|y-{\mathcal{I}}
y\right\|_{0,\hat{I}} \leq \left( \frac{b-a}{2} \right)^{s}\left[ \frac{(k-s)!}{(k+s)!} \right]^{1/2}
\left|y\right|_{s+1,\hat{I}} , \ \ 0 \leq s \leq k.
\end{gather*}
\end{lemma}

We will also utilize the following result, in order to handle the boundary values of 
the (derivative of the) interpolation error.

\begin{lemma}\cite{HSS}\label{prop:HSS}
: Let $\hat{I}=(a,b)$ and assume $y\in H^k(\hat{I})$ for some
integer $k\ge1$. Further, let $\Pi_p y$ be the $L^2(\hat{I})$-projection of $%
y$ onto $\mathbb{P}_p$, $p\ge0$. Then the following error estimates hold,
for $0\le s\le\min\{p+1,k\}$: 
\begin{gather*}
\left\|y-\Pi_p y\right\|_{0,\hat{I}}^2 \le \left( \frac{b-a}{2} \right)^{2s} \frac{(p+1-s)!}{(p+1+s)!}
\left|y\right|_{s,\hat{I}}^2.
\end{gather*}
\begin{gather*}
\left|\left(y-\Pi_py\right)(\pm 1)\right|^2 \le \left( \frac{b-a}{2} \right)^{2s} \frac{1}{2p+1}\frac{(p+1-s)!%
}{(p+1+s)!} \left|y\right|_{s,\hat{I}}^2.
\end{gather*}
\end{lemma}

We close this section noting the following.

\begin{lemma}\cite{ZX}
For the weak derivative of the interpolant, the following holds for $u \in
H^1_0(I)$: 
\begin{equation*}
D_{p-1}({\mathcal{I}} u) = ({\mathcal{I}} u)^{\prime }.
\end{equation*}
Actually $D_{p-1}q = q^{\prime }\ \forall \ q\in\mathbb{P}_p$.
\end{lemma}

\subsection{Error estimates}

The triangle inequality gives
\[
\en{u - u_p}_p \leq \en{u - {\cal{I}} u}_p + \en{{\cal{I}} u - u_p}_p .
\]
The interpolation error will be dealt with in the next section.
In order to deal with the term $\en{{\cal{I}} u - u_p}_p$, we use the following {\emph{error equation}} \cite{ZX}.

\begin{lemma}\cite{ZX}
Let $u$ and $u_p \in \mathcal{M}^\omega_0$ be the solutions of 
(\ref{de})--(\ref{bc}) and (\ref{Auv}) respectively. Then for any
 $v = (v_0,v_b) \in \mathcal{M}^\omega_0$ there holds 
\[
{\mathcal{A}}_p({\mathcal{I}} u - u_p , v ) = {\cal{E}}_1(u,v) + {\cal{E}}_2(u,v) + {\cal{E}}_3(u,v) =: {\cal E}(u,v),
\]
where ${\mathcal{I}}$ is the interpolation operator of Definition \ref{interpolant}, and
\[
{\cal{E}}_1(u,v) :=  \eps_1 \sum_{j=1}^N \langle (u-{\mathcal{I}} u)^{\prime }, (v_0-v_b)\cdot n_j\rangle_{\partial I_j},
\]
\[
{\cal{E}}_2(u,v) := \eps_2 \sum_{j=1}^N \left (u-{\mathcal{I}} u, (b v_0)^{\prime}\right)_{I_j} = \eps_2 \left (u-{\mathcal{I}} u, (b v_0)^{\prime}\right)_I,
\]
\[
{\cal{E}}_3(u,v) := \sum_{j=1}^N (r ({\mathcal{I}} u-u),v_0)_{I_j} = (r ({\mathcal{I}} u-u),v_0)_I.
\]
\end{lemma}

The above result allows us to get a handle on the error $\xi \coloneqq {\cal{I}} u - u_p$, as follows:
\begin{eqnarray*}
\en{\xi}_p^2  &\lesssim& {\cal{A}}_p ({\cal{I}} u - u_p, \xi ) \\
&=& \varepsilon_1 \sum_{j=1}^N \langle (u-{\cal{I}} u)^{\prime }, (\xi_0-\xi_b)\cdot
n_j\rangle_{\partial I_j} + \eps_2  \left(u-{\cal{I}} u, (b\xi_0)^{\prime}\right)_{I} + (r ({\cal{I}} u-u),\xi_0)_I \\
&=&\eps_1 \sum_{j=1}^N \left\{(u-\mathcal{I} u)^{\prime -}_j)
[\xi_0 (x^-_j) - \xi_{b,j}] - (u-\mathcal{I} u)^{\prime +}_{j-1}) [\xi_0
(x^+_{j-1}) - \xi_{b,j-1}] \right\} + \\
&+& \eps_2 \left( u - {\cal{I}} u, (b \xi_0)' \right)_I+ (r(\mathcal{I} u-u), \xi_0 )_I.
\end{eqnarray*}

We note that $\eps_1^{1/2} [\xi_0 (x^+_{j-1}) - \xi_{b,j-1}]$ and $\eps_1^{1/2} [\xi_0 (x^-_j) - \xi_{b,j}]$ may be bounded by $\sigma_j^{-1/2} \en{\xi}_p$.
Moreover,
\[
|  (r(\mathcal{I} u-u), \xi_0 )_I | \leq \Vert \sqrt{r} (\mathcal{I} u-u) \Vert_0 \Vert \xi \Vert_0 \lesssim  \Vert \mathcal{I} u-u  \Vert_0 \en{\xi}_p
\]
hence
\begin{gather} \label{bound_ksi0}
  \begin{split}
    \en{\xi}_p & \lesssim  \varepsilon_1^{1/2} \sum_{j=1}^N
\left(\left|(u-\mathcal{I} u)^{\prime -}_j)\right| + \left|(u-\mathcal{I}
u)^{\prime +}_{j-1})\right| \right) \sigma^{-1/2}_j \\
      & \qquad +  \eps_2 \abs{(\mathcal{I} u-u, (b \xi_0)' )_{I}}
                  + \norm{\mathcal{I} u-u}_0.
  \end{split}
\end{gather}

In order to improve the bound (\ref{bound_ksi0}), we need to choose
our mesh. This choice depends on the regime we are in, so in the next section
we consider each one separately.

\section{Approximation results for the different regimes}
\label{sec:approx}

The design of the mesh (and hence the interpolant) differs from case to case, based on the 
relationship between $\eps_1$, $\eps_2$ (see Table 1). In this Section we consider the three
regimes separately, and provide error estimates that are independent of  $\eps_1, \eps_2$. In particular,
using the triangle inequality we have
\begin{equation}
\label{error0}
\en{u-u_p}_p \leq \en{u - {\cal{I}} u}_p + \en{{\cal{I}} u - u_p}_p = \en{u - {\cal{I}} u}_p + \en{\xi}_p,
\end{equation}
so by the coericity of ${\cal{A}}_p$ we get, recalling that $(u - {\cal{I}} u)$ is continuous,
\begin{align*}
   \en{u - {\cal{I}} u}_p^2
     &\lesssim
        \mathcal{A}_p (u - {\cal{I}} u,u - {\cal{I}} u)  \\
     & \lesssim \eps_1 \norm{D_{p-1}(u - {\cal{I}} u)}_{0,I}^2
                   + \eps_2 (  D^c_{p} (u - {\cal{I}} u), u - {\cal{I}} u)_I
                   +  \norm{ \sqrt{r} (u - {\cal{I}} u) }_{0, I}^2 \\
     & \qquad
         + S(u - {\cal{I}} u,u - {\cal{I}} u )+ S_c(u - {\cal{I}} u,u - {\cal{I}} u) +
                 \abs{u - {\cal{I}} u}_J^2  \\
     & \lesssim
          \eps_1 \norm{(u - {\cal{I}} u)'}_{0,I}^2
        + \eps_2 (  b (u - {\cal{I}} u)', u - {\cal{I}} u)_I
            +  \norm{u - {\cal{I}} u}_{0, I}^2 .
\end{align*}
The middle term above may be absorbed in the last term, after integrating by parts. Hence,
\begin{equation}\label{en_norm_interp_error}
\en{u - {\cal{I}} u}_p  \lesssim  \eps_1^{1/2} \Vert (u - {\cal{I}} u)' \Vert_{0,I}  +  \Vert u - {\cal{I}} u \Vert_{0, I}
\end{equation}

Combining (\ref{bound_ksi0})--(\ref{en_norm_interp_error}), we have
\begin{gather} \label{error1}
  \begin{split}
    \en{u-u_p}_p & \lesssim
       \eps_1^{1/2} \norm{(u - {\cal{I}} u)'}_{0, I}
           + \norm{u - {\cal{I}} u}_{0, I}
           + \eps_2 \abs{ (\mathcal{I} u-u, (b \xi_0)' )_{I}} \\
      & \qquad + \varepsilon_1^{1/2} \sum_{j=1}^N \left(\abs{(u-\mathcal{I} u)^{\prime -}_j)}
               + \abs{(u-\mathcal{I} u)^{\prime +}_{j-1})} \right) \sigma^{-1/2}_j . 
  \end{split}
\end{gather}
The first three terms in (\ref{error1}) will be handled by interpolation estimates, and the last term  by Lemma \ref{prop:HSS},
since the interpolant ${\cal{I}} u $ is constructed as the integral of the Legendre expansion of $u'$ (see \cite{Schwab98}).

\subsection{The regime $\eps_1 \ll \eps_2^2 \ll 1$}

This is the most interesting regime, since the solution will feature two different width boundary layers at the endpoints.
The mesh $\omega$ is chosen
as a variant of the {\emph{Spectral Boundary Layer}} mesh, as defined below.

\begin{definition}\label{SBL_RCD}
Let $\mu _{0},\mu _{1}$ be given by (\ref{mu}). \ For $\kappa >0$, $p\in \mathbb{N}$ 
and $0 < \eps_1 \ll \eps_2^2 \ll 1$, define the Spectral Boundary Layer mesh $\omega _{BL}(\kappa ,p)$ as
\begin{equation*}
\omega _{BL}(\kappa ,p):=\left\{ 
\begin{array}{ll}
\omega =\{0,1\} & \text{if }\kappa p\varepsilon_1 \geq 1/2 \\ 
\omega =\{0,\kappa p\mu _{0}^{-1},1-\kappa p\mu _{1}^{-1},1\} & \text{if }\kappa p\varepsilon_2 <1/2
\end{array}
\right. .
\end{equation*}
\end{definition}

The following interpolation error estimates were shown in  \cite{SXArXiv}:
\begin{equation}
\left. 
\begin{array}{r@{\ }l}
  \norm{\left( u-{\cal{I}} u \right)'}_{0, I}
   & \lesssim \eps_1^{-1/2}e^{-\gamma _{1}p} \\ 
  \norm{u-{\cal{I}} u}_{0, I}
   & \lesssim e^{-\gamma _{1}p} \\ 
  \norm{u-{\cal{I}} u}_{0, (0,\mu_0^{-1})}
   & \lesssim \eps_1^{1/2}e^{-\gamma _{1}p} \\
  \norm{u-{\cal{I}} u}_{0, (1-\mu_1^{-1},1)}
   & \lesssim \eps_1^{1/2} \eps_2^{-1}  e^{-\gamma _{1}p}
\end{array}%
\right\}.  \label{interp1d_c}
\end{equation}
for some constant $\gamma _{1}>0$ independent of $\eps_1, \eps_2$.

The main result of this section is the following.

\begin{proposition}
\label{prop:regime3}
Let $u, u_p$ be the solutions of (\ref{de})--(\ref{bc}) and (\ref{Auv}), respectively and assume $\eps_1 \ll \eps_2^2$. Then
there exists a positive constant $\sigma$ such that
\[
\en{u-u_p}_p \lesssim e^{-\sigma p}.
\]
\end{proposition}
\begin{proof}
The proof is separated in two cases.

\paragraph{\textit{Case 1}: $\kappa p \eps_1 \geq 1/2$.}
In this case the mesh consists of one element and there holds $\eps_1^{-1} \leq 2 \kappa p$. The solution satisfies
(\ref{Classical3}) and we have from (\ref{error1}), with the aid of  (\ref{interp1d_c}),
\[
\en{u-u_p}_p \lesssim e^{-\gamma_1 p}  + \eps_2 \vert (\mathcal{I} u-u, (b \xi_0)' )_{I} \vert   + \varepsilon_1^{1/2} \sum_{j=1}^N
\left(\left|(u-\mathcal{I} u)^{\prime -}_j)\right| + \left|(u-\mathcal{I}
u)^{\prime +}_{j-1})\right| \right) \sigma^{-1/2}_j .
\]
The second term above satisfies
\[
\eps_2 \vert (\mathcal{I} u-u, (b \xi_0)' )_{I} \vert \lesssim  \eps_2 \norm{\mathcal{I} u-u}_{0,I}
    \norm{\xi'_0}_{0,I} 
\lesssim \eps_2  \norm{\mathcal{I} u-u}_{0,I} \frac{ p^2 }{|I|} \norm{\xi_0}_{0,I} \lesssim e^{-\gamma_1 p},
\]
where (\ref{interp1d_c}) and an inverse inequality (\cite[Thm. 3.91]{Schwab98}) were used.

It remains to consider the term
\begin{eqnarray*}
 \varepsilon_1^{1/2} \sum_{j=1}^N
\left(\left|(u-\mathcal{I} u)^{\prime -}_j)\right| + \left|(u-\mathcal{I} u)^{\prime +}_{j-1})\right| \right) \sigma^{-1/2}_j 
&=&  \varepsilon_1^{1/2} \left(\left|(u-\mathcal{I} u)^{\prime -}_1)\right| + \left|(u-\mathcal{I}
u)^{\prime +}_0)\right| \right) \sigma^{-1/2}_1 \\
&\leq&  \frac{1}{2p+1}\frac{(p+1-s)!}{(p+1+s)!} \left| u' \right|_{s,{I}}^2,
\end{eqnarray*}
by Lemma \ref{prop:HSS}. Choosing $s = \lambda (p+1)$ with $\lambda \in (0,1)$ arbitrary, we have
\[
 \frac{1}{2p+1}\frac{(p+1-s)!}{(p+1+s)!} \left| u' \right|_{s,{I}}^2 \leq \frac{ ((p+1)-\lambda (p+1))!}{((p+1)+\lambda (p+1))!}
 \left| u \right|_{\lambda (p+1)+1,{I}}^2
\]
with
\[
 \left| u \right|_{\lambda (p+1)+1,I}^2 = \int_{I} [u^{(\lambda (p+1)+1}(x)]^2 dx \lesssim
 K_1^{2\lambda (p+1)+1} \eps_1^{-2\lambda (p+1)+1} \lesssim K_1^{2\lambda (p+1)+1}  (2 \kappa p)^{2 \lambda (p+1)+1}.
\]
Hence,
\[
 \frac{1}{2p+1}\frac{(p+1-\lambda (p+1))!}{(p+1+\lambda (p+1))!} \left| u \right|_{\lambda (p+1)+1,{I}}^2 \lesssim
\frac{(p+1-\lambda (p+1))!}{(p+1+\lambda (p+1))!} K_1^{2\lambda (p+1)+1}  (2 \kappa p)^{2\lambda (p+1)+1}.
\]
The factorial term simplifies, with the aid of Stirling's formula (see, e.g., \cite{MXO}),
and we obtain
\begin{align*}
  & \frac{(p+1-\lambda (p+1))!}{(p+1+\lambda (p+1))!} \left| u \right|_{\lambda (p+1)+1,{I}}^2 \\
  & \qquad \leq \left[ \frac{(1-\lambda)^{1-\lambda}}{(1+\lambda)^{1+\lambda}}\right]^{p+1}
              (p+1)^{-2 \lambda (p+1)+1} e^{2 \lambda (p+1)}
              K_1^{2 \lambda (p+1)}  (2 \kappa p)^{2\lambda (p+1)} \\
  & \qquad \lesssim \left[ \frac{(1-\lambda)^{1-\lambda}}{(1+\lambda)^{1+\lambda}}
         (2 \kappa  K_1 e)^{2 \lambda} \right]^{p+1} \lesssim e^{-\beta p},
\end{align*}
where $\beta = |\ln q |$, $q = \frac{(1-\lambda)^{1-\lambda}}{(1+\lambda)^{1+\lambda}} \in (0, 1)$,
and we selected the constant $\kappa$ in the definition of the mesh as 
$\kappa = \frac{1}{2 K_1 e}$ .
The desired result, in this case, follows.

\paragraph{\textit{Case 2}: $\kappa p \eps_2 < 1/2$.}
We recall the decomposition (\ref{decomposition})
\[
u=u_{S}+u^{\pm}_{BL}+u_{R}, 
\]
with the bounds (\ref{Classical3})--(\ref{rM1bound}) being true. The mesh consists of the three elements
\begin{equation}
\label{intervals}
I_1 = (0,\kappa p \mu_0^{-1}), I_2 = (\kappa p \mu_0^{-1}, 1 - \kappa p \mu_1^{-1}), I_3 = (1-\kappa p \mu_1^{-1}, 1),
\end{equation}
and each component in the decompostion is approxiated separately. In particular, the approximation of the smooth part $u_S$ is achieved by
a degree $p$ interpolant on the entire interval, just like in {\textit{Case 1}} above -- the details are omitted. The boundary layer $u_{BL}^{-}$ on 
the left will be approximated by a degree $p$ interpolant on $I_1$ and by its linear interpolant on $I_2 \cup I_3$. Similarly, the boundary layer 
$u_{BL}^{+}$ on the right is approximated by its linear interpolant on $I_1 \cup I_2$ and by a degree $p$ interpolant on $I_3$. The remainder
(cf. (\ref{rM1bound})) will not be approximated, since it is already exponentially small. Below we provide the details only for $u_{BL}^-$ since the
approximation of $u_{BL}^+$ is very similar.

So, in $I_2 \cup I_3$ we have for $u_{BL}^{-}$
\begin{align*}
 & \en{u_{BL}^{-} - {\cal{I}}_1 u_{BL}^{-} }_{p,I_2 \cup I_3}
    \lesssim \en{u^-_{BL}}_{p,I_2 \cup I_3} + \en{ {\cal{I}}_1 u_{BL}^{-}}_p
    \lesssim  \en{ {\cal{I}}_1 u_{BL}^{-}}_p \\
 & \qquad \lesssim \eps_1^{1/2} \norm{( u_{BL}^-)'}_{0,I_2 \cup I_3} + \norm{\mathcal{I}_1 u_{BL}^{-}}_{0,I_2 \cup I_3}
    \lesssim \eps_1^{1/2} \eps_2^{-1} e^{- \frac{\kappa p}{\eps_2 \mu_0}}
    \lesssim e^{-\beta p},
\end{align*}
where we used (\ref{p_norm}), (\ref{uBL1abound}) and (\ref{mu_a}).

Next, we consider $u_{BL}^{-}$ on $I_1$, and we have by (\ref{error1})--(\ref{intervals})
\begin{gather}\label{summ}
  \begin{split}
    \en{u_{BL}^{-} - {\cal{I}}_p u_{BL}^{-}}_{p,I_1} & \lesssim 
      e^{-\gamma_1 p} + \eps_2 \abs{\left(\mathcal{I} u_{BL}^{-} -u_{BL}^{-} , (b \xi_0)' \right)_{I_1}} \\
      & \qquad + \varepsilon_1^{1/2} 
         \left(\abs{(u_{BL}^{-} -\mathcal{I} u_{BL}^{-} )^{\prime -}_1)}
                 + \abs{(u_{BL}^{-} -\mathcal{I}
u_{BL}^{-} )^{\prime +}_{0})} \right) \sigma^{-1/2}_1.
  \end{split}
\end{gather}
The second term in (\ref{summ}) satisfies
\[
\eps_2 \vert (\mathcal{I} u_{BL}^{-} -u_{BL}^{-} , (b \xi_0)' )_{I_1} \vert \lesssim 
\eps_2 \Vert {\cal{I}} u_{BL}^{-}  - u_{BL}^{-}  \Vert_{0,I_1}  \Vert \xi'_0 \Vert_{0,I_1} 
\lesssim \eps_2 \Vert {\cal{I}} u_{BL}^{-}  - u_{BL}^{-}  \Vert_{0,I_1} \frac{ p^2 }{|I_1|} \Vert\xi_0 \Vert_{0,I_1},
\]
and since $\vert I_1 \vert=\kappa p \mu_0^{-1}$, we have from  (\ref{mu_a}) that $\frac{\eps_2}{\vert I_1 \vert} \leq \eps_2 \mu_0 \lesssim 1$.
Hence, by (\ref{interp1d_c})
\begin{equation} \label{summ1}
\eps_2 \vert (\mathcal{I} u_{BL}^{-} -u_{BL}^{-} , (b \xi_0)' )_{I_1} \vert \lesssim e^{-\beta p}.
\end{equation}
For the last term in (\ref{summ}), we first note that
$\eps_1^{1/2} \sigma_1^{-1/2} = \kappa / \mu_0$, so by Lemma \ref{prop:HSS}
\begin{align*}
  & \eps_1^{1/2} \left(\abs{(u_{BL}^{-} -\mathcal{I} u_{BL}^{-} )^{\prime -}_1)}
       + \abs{(u_{BL}^{-} -\mathcal{I} u_{BL}^{-} )^{\prime +}_{0})} \right)
      \sigma^{-1/2}_1 \\
  & \qquad
    \leq  \frac{\kappa}{\mu_0}\left( |I_1|/2 \right)^{2s} \frac{1}{2p+1}\frac{(p+1-s)!}{(p+1+s)!} \left| (u_{BL}^{-})' \right|_{s,{I_1}}^2.
\end{align*}
 Choosing $s = \lambda_1 (p+1)$ with $\lambda_1 \in (0,1)$ arbitrary, we have
\begin{align*}
 & \frac{1}{2p+1}\frac{(p+1-s)!}{(p+1+s)!} \abs{(u_{BL}^{-})'}_{s,{I_1}}^2
     \leq \frac{ ((p+1)-\lambda_1 (p+1))!}{((p+1)+\lambda_1 (p+1))!} 
 \left| u_{BL}^{-} \right|_{ \lambda_1 (p+1)+1,{I_1}}^2 \\
  & \qquad
    \lesssim \left[ \frac{(1-\lambda_1)^{1-\lambda_1}}{(1+\lambda_1)^{1+\lambda_1}}\right]^{p+1}
(p+1)^{-2 \lambda_1 (p+1)+1} e^{2 \lambda_1 (p+1)+1}  \left| (u_{BL}^{-}) \right|_{ \lambda_1 (p+1)+1,{I_1}}^2,
\end{align*}
and since
\[
 \abs{u_{BL}^{-}}_{ \lambda_1 (p+1)+1,I_1}^2
    = \int_{I_1} \bigl[(u_{BL}^{-})^{ \lambda_1 (p+1)+1}(x)\bigr]^2 dx
    \lesssim \tilde{K}^{2( \lambda_1 (p+1)+1)} \eps_2^{-2( \lambda_1 (p+1)+1)},
\]
we arrive at
\begin{align*}
 & \eps_1^{1/2} \left(\abs{(u_{BL}^{-} -\mathcal{I} u_{BL}^{-} )^{\prime -}_1)}
      + \abs{(u_{BL}^{-} -\mathcal{I} u_{BL}^{-} )^{\prime +}_{0})} \right) \sigma^{-1/2}_1  \\
 & \qquad
    \lesssim  \frac{\kappa}{\mu_0}\left(\frac{\kappa p \mu_0^{-1}}{2} \right)^{2 \lambda_1(p+1)}
       \left[ \frac{(1-\lambda_1)^{1-\lambda_1}}{(1+\lambda_1)^{1+\lambda_1}}\right]^{p+1} \\
 & \qquad\qquad\qquad\qquad
         \times (p+1)^{-2 \lambda_1 (p+1)} e^{2 \lambda_1 (p+1)+1}  \tilde{K}^{2(\lambda_1 (p+1)+1)}
      \eps_2^{-2( \lambda_1 (p+1)+1)} \\
 & \qquad
    \lesssim \left( \mu_0 \eps_2\right)^{-2\lambda_1 (p+1)+1} 
       \left[ \frac{(1-\lambda_1)^{1-\lambda_1}}{(1+\lambda_1)^{1+\lambda_1}} \left( \frac{\kappa e \tilde{K}}{2} \right)^{2 \lambda_1} \right]^{p+1}
 \lesssim e^{-\tilde{\beta} p},
\end{align*}
where $\tilde{\beta} = |\ln \tilde{q} |, \tilde{q} = \frac{(1-\lambda_1)^{1-\lambda_1}}{(1+\lambda_1)^{1+\lambda_1}} \in (0, 1)$. Here, we selected the 
constant $\kappa$ in the definition of the mesh as $\kappa = \frac{2}{\tilde{K} e}$, and made use of (\ref{mu_a}).

Combining the above gives the desired result for $u_{BL}^-$.
\end{proof}

\subsection{The regime $\eps_1 \gg \eps_2^2$}

This is the typical reaction-diffusion SPP, with layers of width $O(1/\sqrt{\eps_1})$ near both endpoints,
and the {\emph{Spectral Boundary Layer}} mesh \cite{MXO} is defined below.

\begin{definition}\label{SBL_RD}
For $\kappa >0$, $p\in \mathbb{N}$ 
and $0<\eps_{2}^2 \ll \eps_{1} < 1$, define the Spectral Boundary Layer mesh $\omega _{BL}(\kappa ,p)$ as
\begin{equation*}
\omega _{BL}(\kappa ,p):=\left\{ 
\begin{array}{ll}
\omega =\{0,1\} & \text{if }\kappa p\varepsilon_1 \geq 1/2 \\ 
\omega =\{0,\kappa p\eps^{1/2}_1,1-\kappa p\eps _{1}^{1/2},1\} & \text{if }\kappa p\varepsilon_2 <1/2
\end{array}
\right. .
\end{equation*}
\end{definition}

Moreover,  the following error estimates were shown in  \cite{melenk97}:
\begin{equation}
\left. 
\begin{array}{r@{\ }l}
  \norm{\left( u-{\cal{I}} u \right) ^{\prime }}_{0, I}
    & \lesssim \eps_1^{-1/2}e^{-\gamma _{2}p} \\ 
  \norm{u-{\cal{I}} u}_{0,I}
    & \lesssim e^{-\gamma _{2}p} \\ 
  \norm{u-{\cal{I}} u}_{0, (I_1 \cup I_3)}
    & \lesssim \eps_1^{1/2}e^{-\gamma _{2}p}
\end{array}%
\right\}.  \label{interp1d_b}
\end{equation}
for some constant $\gamma _{2}>0$ independent of $\eps_1, \eps_2$.

The following is the main result of this subsection.
\begin{proposition}
\label{prop:regime1}
Let $u, u_p$ be the solutions of (\ref{de})--(\ref{bc}) and (\ref{Auv}), respectively and assume $\eps_1 \gg \eps_2^2$. Then
there exists a positive constant $\sigma$ such that
\[
\en{u-u_p}_p \lesssim e^{-\sigma p}.
\]
\end{proposition}
%

\subsection{The regime $\eps_1 \ll \eps_2=1$}

This is the typical convection-diffusion SPP, with a layer of width $O(1/\eps_1)$ near the right endpoint. The
{\emph{Spectral Boundary Layer}} mesh  is defined accordingly.

\begin{definition}\label{SBL_CD}
For $\kappa >0$, $p\in \mathbb{N}$ 
and $0<\eps_{1} < 1$, define the Spectral Boundary Layer mesh $\omega _{BL}(\kappa ,p)$ as
\begin{equation*}
\omega _{BL}(\kappa ,p):=\left\{ 
\begin{array}{ll}
\omega =\{0,1\} & \text{if }\kappa p\eps_1 \geq 1/2 \\ 
\omega =\{0,1-\kappa p\eps _{1},1\} & \text{if }\kappa p\eps_1 \leq 1/2
\end{array}
\right. .
\end{equation*}
\end{definition}

The following interpolation error estimates were shown in  \cite[Lemma 2.4]{MS}:
\begin{equation}
\left. 
\begin{array}{r@{\ }l}
  \norm{\left( u-{\cal{I}} u \right)'}_{0, I}
    & \lesssim \eps_1^{-1/2}e^{-\gamma _{3}p} \\ 
  \norm{u-{\cal{I}} u}_{0, I}
    & \lesssim e^{-\gamma _{3}p} \\ 
  \norm{u-{\cal{I}} u}_{0,(1-\eps_1,1)}
    & \lesssim \eps_1^{1/2}e^{-\gamma _{3}p}
\end{array}%
\right\}.  \label{interp1d}
\end{equation}
for some constant $\gamma _{3}>0$ independent of $\eps_1$.

The main result of this section is given below.

\begin{proposition}
\label{prop:regime2}
Let $u, u_p$ be the solutions of (\ref{de})--(\ref{bc}) and (\ref{Auv}), respectively and assume $\eps_1 \ll \eps_2=1$. Then
there exists a positive constant $\sigma$ such that
\[
\en{u-u_p}_p \lesssim e^{-\sigma p}.
\]
\end{proposition}
%

%
%
%
%
%
%
%

\section{Numerical Results}

\label{NR} In this section we show the results of numerical computaitons in order to illustrate the theory. We focus on the 
reaction-convection-diffusion case, and consider the model problem
\begin{gather}
  -\eps_1 u'' + \eps_2 \cos(x) u' + (1+x) u = e^{x}\,, \ \ u(0)=u(1)=0.
\end{gather}
We use the 
{\emph{Spectral Boundary Layer}} mesh of Definition \ref{SBL_RCD},
with the  constant $\kappa=1$ in the definition of the mesh (other choices gave similar results). 

In Figure \ref{soln} below, we show the WG solution computed with $p=4, \eps_1 = 10^{-5}, \eps_2 = 10^{-2}$.
The dots correspond to the values at the nodes, and the discontinuous nature is visible.

\begin{figure}[h]
\begin{center}
\includegraphics[width=0.65\textwidth]{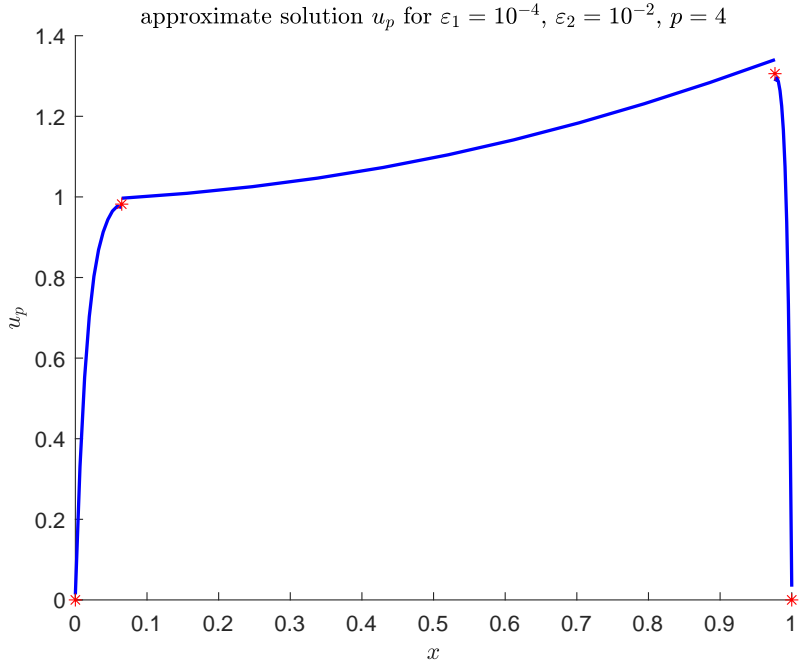}
\end{center}
\caption{The approximate solution.}
\label{soln}
\end{figure}

Using twice as high polynomial degree, we obtain a reference solution (since there is no exact one available) and we calculate the (estimated) error in 
the approximation for the $\en{\cdot}_p$ (energy) norm.

\begin{figure}[h]
\begin{center}
\includegraphics[width=0.65\textwidth]{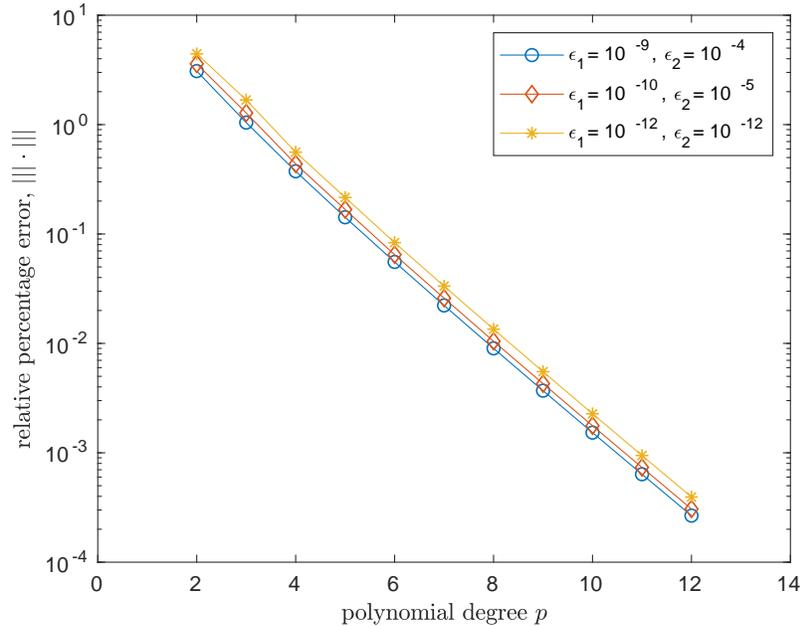}
\end{center}
\caption{Convergence of the method -- all regimes.}
\label{f1}
\end{figure}

\begin{figure}[h]
\begin{center}
\includegraphics[width=0.65\textwidth]{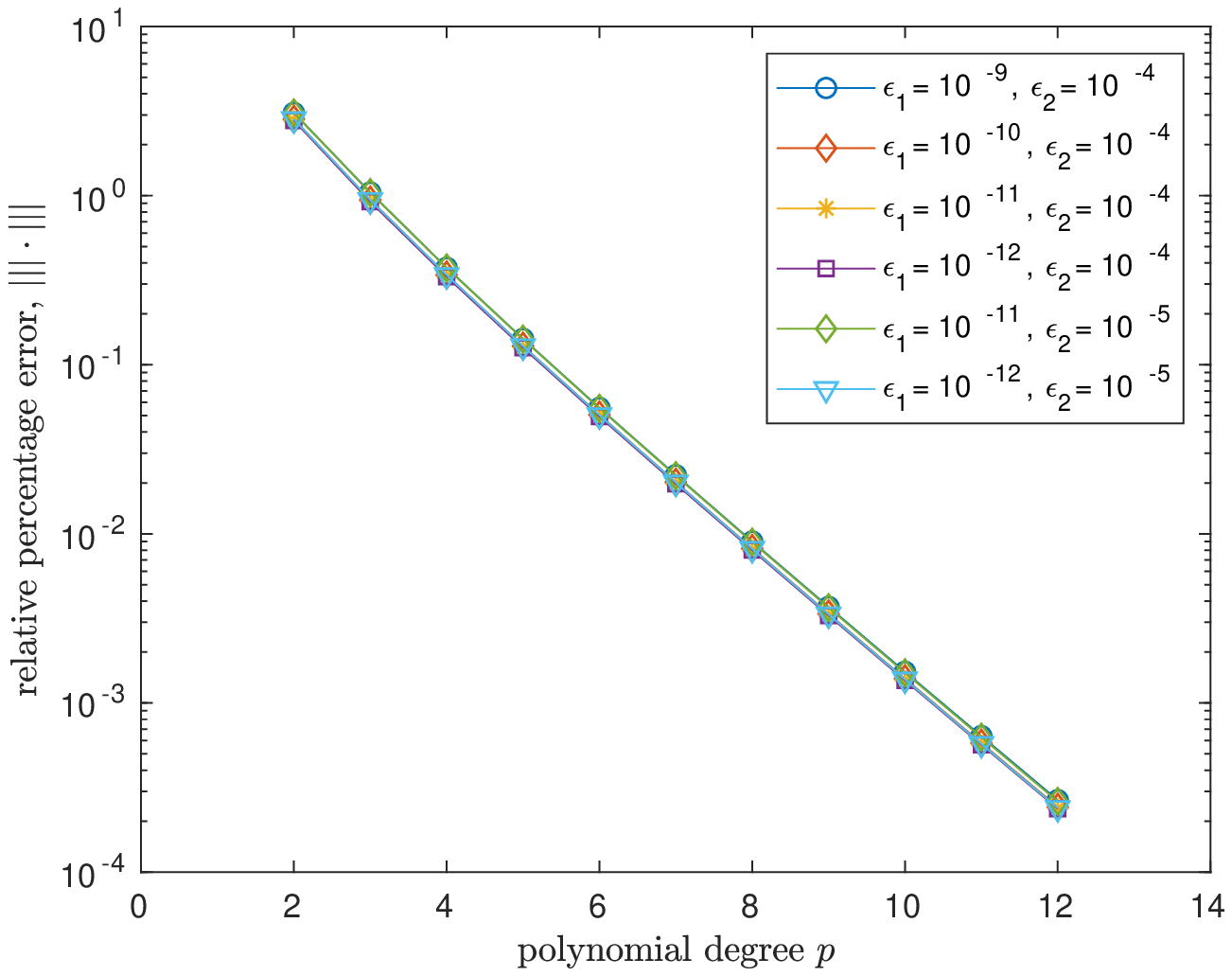}
\end{center}
\caption{Convergence of the method -- regime 1.}
\label{f2}
\end{figure}

\begin{figure}[h]
\begin{center}
\includegraphics[width=0.65\textwidth]{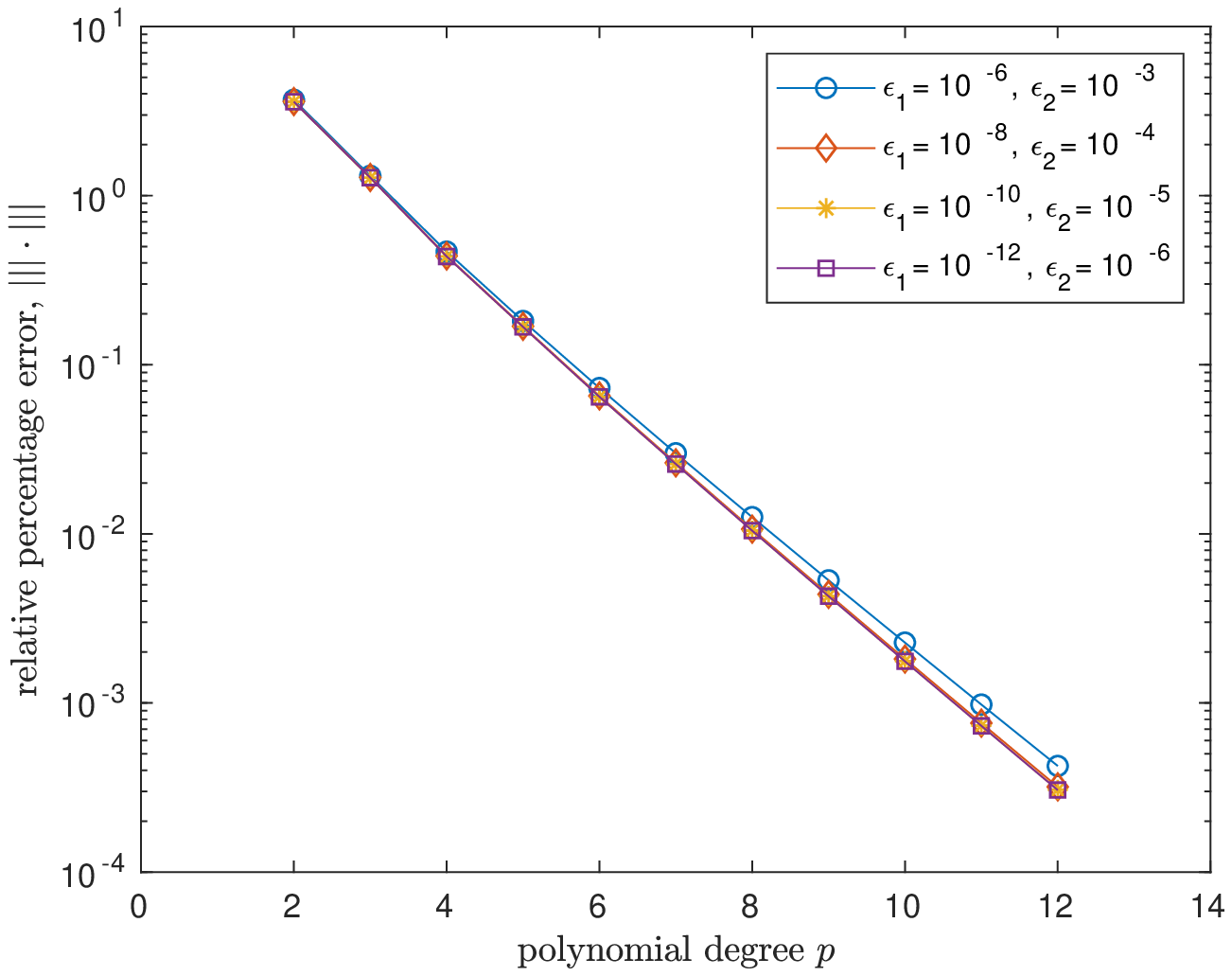}
\end{center}
\caption{Convergence of the method -- regime 2.}
\label{f3}
\end{figure}

\begin{figure}[h]
\begin{center}
\includegraphics[width=0.65\textwidth]{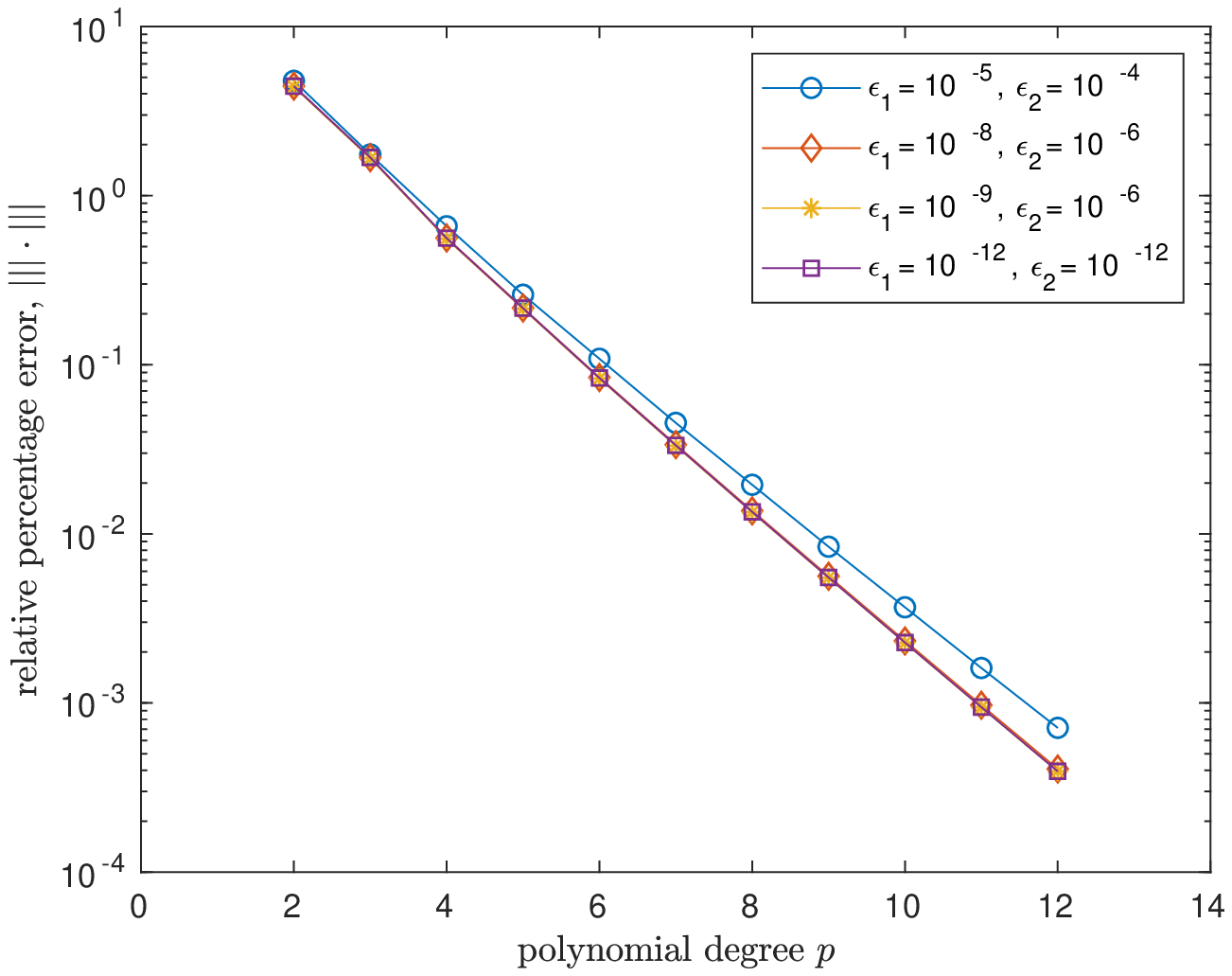}
\end{center}
\caption{Convergence of the method -- regime 3.}
\label{f4}
\end{figure}

In Figures \ref{f1} -- \ref{f4}, we show semi-log plots of the percentage relative errors (in the energy norm), versus the
polynomial degree, for various combinations of $\eps_1, \eps_2$ which represent all three regimes.

We observe in {\emph{all}} cases, that the method converges uniformly and at an exponential rate, as the theory predicts.



%
%



\end{document}